\documentclass[11pt]{amsart}
\usepackage{amsfonts}
\usepackage{amsmath}
\usepackage{amssymb}
\usepackage{setspace}
\usepackage[dvips]{graphicx}
\newtheorem{theorem}{Theorem}
\newtheorem{prop}[theorem]{Proposition}

\newcommand{\leb}{\mathcal{L}}
\newcommand{\supp}{\textrm{Supp}}

\newcommand{\e}{\varepsilon}

\newcommand{\Q}{\mathcal{Q}}

\newcommand{\R}{\mathcal{R}}
\newcommand{\D}{\mathcal{D}}
\newcommand{\rot}{\mathbf{R}}
\newcommand{\RR}{\mathbb{R}}

\title{Resonance between Cantor sets}
\author{Yuval Peres}\thanks{Research of Y. Peres supported in part by NSF grant
   DMS-0605166}

\address{Yuval Peres,
Microsoft Research, Redmond and Departments of Statistics and
Mathematics, University of California, Berkeley.}
\email{peres@stat.Berkeley.edu}

\author{Pablo Shmerkin}
\thanks{Research of P. Shmerkin supported in part by the Academy of Finland and Microsoft Research}
\address{Pablo Shmerkin, Departments of Mathematics and  Statistics, University of Jyv\"askyl\"a.}
\email{shmerkin@maths.jyu.fi}

\begin{document}

\begin{abstract}
Let $C_a$ be the central Cantor set obtained by removing a central
interval of length $1-2a$ from the unit interval, and continuing
this process inductively on each of the remaining two intervals.
We prove that if $\log b/\log a$ is irrational, then
\[
\dim(C_a+C_b) = \min(\dim(C_a) + \dim(C_b),1),
\]
where $\dim$ is Hausdorff dimension. More generally, given two
self-similar sets $K,K'$ in $\RR$ and a scaling parameter $s>0$, if the dimension of the
arithmetic sum $K+sK'$ is strictly smaller than $\dim(K)+\dim(K')
\le 1$ (``geometric resonance''), then there exists $r<1$ such
that all contraction ratios of the similitudes defining $K$ and
$K'$ are powers of $r$ (``algebraic resonance''). Our method also
yields a new result on the projections of planar self-similar sets
generated by an iterated function system that includes a scaled
irrational rotation.
\end{abstract}

\maketitle

\section{Introduction and statement of results}

In many situations in dynamical systems and geometric measure
theory one encounters results which are valid for almost every
member of a parametrized family. In general, those results are
sharp in the sense that the set of exceptions may have positive,
or even full, dimension. But when the constructions are done in a
dynamical or geometrically regular way, there is often some
countable set of parameters, usually arising from an algebraic
condition, where the result fails, and it is natural to conjecture
that these are all parameters for which the result fails, see for
example \cite{furstenberg}, \cite{keane-smorodinsky-solomyak},
\cite{pollicott-simon}, \cite{sixtyyears},
\cite{keane-simon-solomyak}. However, there are very few cases
where the set of exceptions has been explicitly determined. In
this paper we will find the exact set of parameters for which
certain kind of resonance between two self-similar sets occurs.
This set is countable and given by a simple, natural algebraic
relation.

Let $C_a$ be the central Cantor set obtained by removing a middle
interval of length $1-2a$ from the unit interval, and continuing
this process inductively on each of the remaining intervals.
 Explicitly,
\[
C_a = \left\{  (1-a)\sum_{i=0}^\infty \omega_i a^i :
\omega_i\in\{0,1\}\right\}.
\]
It is well known that $\dim(C_a) = \log 2/\log(1/a)$.

Our main result is the following:

\begin{theorem} \label{th:sumcantor}
If $\log b/\log a$ is irrational, then
\begin{equation} \label{eq:sumcantor}
\dim(C_a+C_b) = \min(\dim(C_a)+\dim(C_b),1).
\end{equation}
\end{theorem}

When $\log b/\log a$ is rational and $\dim(C_a)+\dim(C_b)\le 1$,
the equality (\ref{th:sumcantor}) does not hold.  This is well
known; we sketch a proof after (\ref{drop}) below.

An {\bf iterated function system} (or i.f.s. for short) on a complete metric
space $\mathbb{X}$ is a finite
collection of maps $\{ f_1,\ldots, f_n\}$ from $\mathbb{X}$ to itself, such that all $f_j$ are contractions
(have  Lipschitz constant strictly less than $1$.) The \textbf{attractor} of this i.f.s.\ is the unique
nonempty compact set $E$ such that
\[
E = \bigcup_{i=1}^n f_i(E).
\]
In particular, $C_a$ is the attractor of the
i.f.s. $\{ a x, \, \, a x + (1-a)\}$ in $\RR$.
When $\mathbb{X}=\mathbb{R}^n$ and the maps are similitudes, we say
that $E$ is a \textbf{self-similar} set. See \cite{falconer2} for
background on iterated function systems and self-similar sets.

Theorem \ref{th:sumcantor} generalizes to arbitrary self-similar
sets in $\mathbb{R}$:

\begin{theorem} \label{th:sumselfsimilar}
Let $K$ (resp. $K'$) be the attractor of the i.f.s. $\{ r_i x +
t_i \}_{i=1}^n$ (resp. $\{r'_i x + t'_i\}_{i=1}^{n'}$). If there
exist $j,j'$ such that $\log(|r_j|)/\log(|r'_{j'}|)$ is
irrational, then
\[
\dim(K+K') = \min(\dim(K)+\dim(K'),1).
\]
\end{theorem}

Next, we define two notions  of {\em resonance}, and restate Theorem \ref{th:sumselfsimilar}
using these notions. Let
$\{ r_i x+t_i\}_{i=1}^n$ and $\{ r'_i
x+t'_i\}_{i=1}^{n'}$, be two i.f.s. consisting of similitudes, with attractors $K$ and $K'$ respectively.
We say that these  two i.f.s. are \textbf{algebraically resonant} if for all
 $j \le n$ and $j' \le n'$, the ratio $\log|r_j|/\log|r'_{j'}|$ is rational.
 We say that the two i.f.s. are \textbf{geometrically resonant}
  if there is some $s>0$ such that
\[
\dim(K+sK') < \dim(K)+\dim(K') \le 1 \,.
\]

Theorem \ref{th:sumselfsimilar} implies that two geometrically
resonant iterated function systems  must also be algebraically
resonant. Below we will state a converse under a mild separation
condition; see Theorem \ref{th:exceptions}.

The paper is structured as follows. In the remainder of this
section we discuss some background related to Theorem
\ref{th:sumcantor}. Section \ref{sec:preliminaries} contains a discrete version of the Marstrand
Projection Theorem,
which will be needed for the main proofs. The
proof of Theorem \ref{th:sumcantor} is in Section
\ref{sec:proofmaintheorem}, while Section
\ref{sec:remainingproofs} contains the proofs of the additional
results stated at the end of this section. We finish the paper
with some remarks and open questions in Section \ref{sec:remarks}.

\subsection{Background}

The study of the arithmetic sums of Cantor sets is a classical
topic in real analysis, motivated in part by its rich applications
to Diophantine approximation, the study of homoclinic bifurcations
in dynamical systems, and other topics. The basic problem is to
compute or estimate the size of the arithmetic sum $K+K'$ in terms
of the sizes of $K$ and $K'$. ``Size'' can mean different things,
but in this paper we focus on Hausdorff dimension.

Let $K, K'$ be two compact subsets of the real line. Observe that
\[
K + K' \mbox{ \rm is congruent to } \sqrt{2}P_{\pi/4}(K\times K'),
\]
where $P_\theta$ denotes the orthogonal projection onto the  line $\{(t\cos \theta, t\sin\theta)\}_{t \in \RR}$
making angle $\theta$ with the $x$-axis. Let
$
\gamma = \dim(K) + \dim(K').
$
By Marstrand's projection theorem (see \cite{mattila}, Chapter
9) we have
\begin{eqnarray*}
\gamma \le 1 & \Longrightarrow & \dim(P_\theta(K\times K')) =
\gamma \,\textrm{ for almost
every } \theta\in [0,\pi),\\
\gamma > 1 & \Longrightarrow & \leb(P_\theta(K\times K'))>0
\,\textrm{ for almost every } \theta\in [0,\pi),
\end{eqnarray*}

where $\leb$ denotes Lebesgue measure. Unfortunately, Marstrand's
theorem gives no information about specific values of $\theta$,
and therefore cannot be directly applied to obtain information
about sums of Cantor sets. However, it does  support the
heuristic principle that ``typically'',
\[
\dim(K+K') = \min(\dim(K)+\dim(K'),1).
\]
Note that our definition of geometrical resonance for the iterated function systems that generate $K$ and $K'$ (in the case
where the sum of the dimensions is at most $1$) requires
at least one exception to the projection theorem for the set
$K\times K'$.
% Variants of Marstrand's theorem have been applied to justify this heuristic principle in some cases.
We now discuss some earlier results that
motivated our work.

\begin{figure} \label{fig:dimdrop}
\begin{center}
\includegraphics[width=0.75\textwidth]{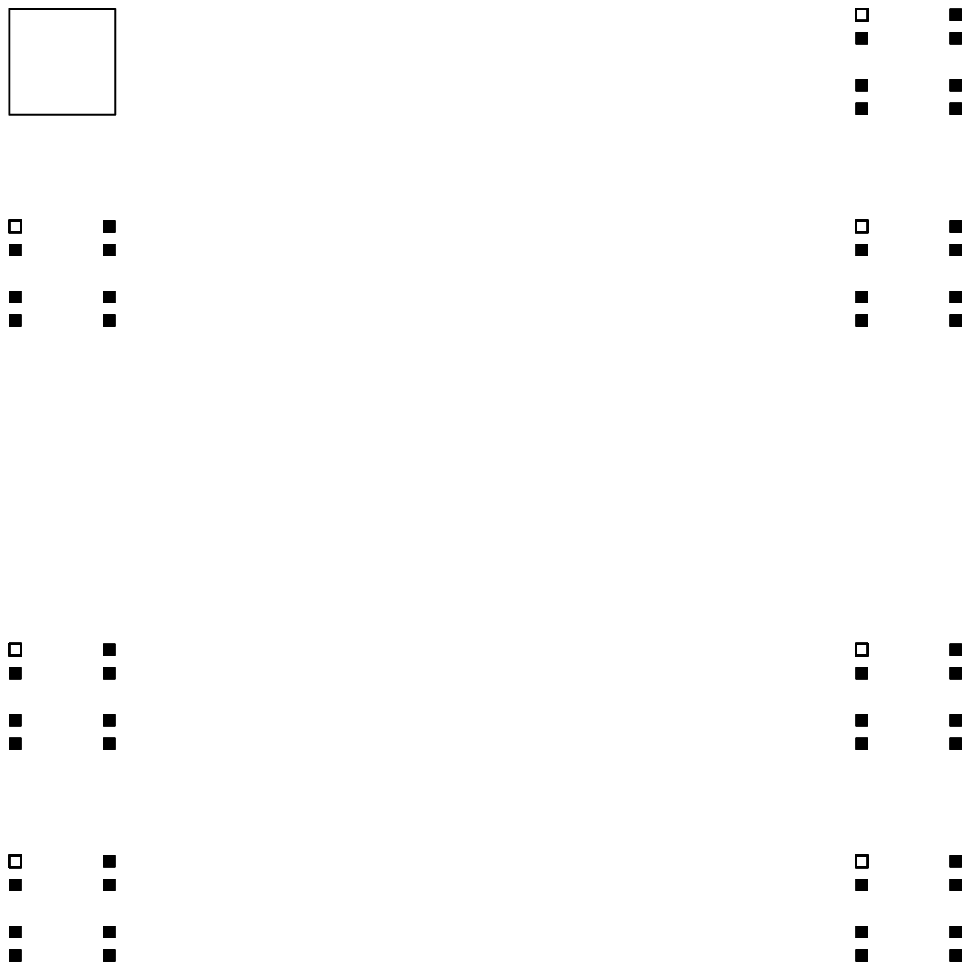}
\end{center}
\textbf{Figure 1}. The structure of $C_{1/9}\times C_{1/3}$: for
each white square there is a black square of the same level, and
with the same projection onto the line $\{
(t,t)\}_{t\in\mathbb{R}}$. This results in a dimension drop for
the projection, which is affinely equivalent to the arithmetic sum
$C_{1/9}+C_{1/3}$.
\end{figure}

The following result is due to Peres and Solomyak
(\cite{peres-solomyak-cantor}, Theorem 1.1), generalizing an
earlier result of Solomyak  \cite{solomyak-indagationes}:
\begin{theorem} \label{th:peressolomyak} (Peres and Solomyak)
Let
\[
\gamma(a) = \dim(C_a) = \log(2)/\log(1/a).
\]
Given a fixed compact set $K\subset\mathbb{R}$, the following holds
for almost every $a\in (0,1/2)$:
\begin{eqnarray}
\gamma(a)+\dim(K) \le 1 & \Longrightarrow & \dim(C_a+ K) =
\gamma(a)+\dim(K),
\label{eq:sumcantorcompactdim}\\
\gamma(a)+\dim(K) > 1 & \Longrightarrow & \leb(C_a+ K)
>0.
\end{eqnarray}
\end{theorem}
In this theorem, the Cantor sets $C_a$ can be replaced by more
general homogeneous self-similar sets  (see
\cite{peres-solomyak-cantor} for details). If the ratio of $\log
a$ and $\log b$ is rational, then
\begin{equation} \label{drop}
\dim(C_a+C_b) < \dim(C_a) + \dim(C_b),
\end{equation}
This is folklore and can be seen as follows: if $\log a/\log b$ is
rational, then there are $0<r<1$ and positive integers $m, n$ such
that $a=r^m,\, b=r^n$. By iterating $n$ times the construction of
$C_a$ and $m$ times the construction of $C_b$, we can decompose
the product $C_a\times C_b$ into a union of cylinder sets whose
convex hulls are squares. Moreover, two of these squares are
$[0,r^{mn}]\times [1-r^{mn},1]$ and $[1-r^{mn},1]\times
[0,r^{mn}]$.
%(corresponding to the first and last intervals in
%the constructions of $C_a$ and $C_b$, respectively).
The projections of these squares on the line $\{(t,t)\}_{t \in \RR}$ coincide. Thus, if  we
delete one of them from the construction of $C_a\times
C_b$, and delete the corresponding
square at every level, we obtain a   self similar set $K \subset C_a \times C_b$  with
$P_{\pi/4}(K)=P_{\pi/4}(C_a\times C_b)$ and $\dim(P_{\pi/4}(K)) \le \dim(K) <\dim(C_a \times C_b)$.
See Figure \ref{fig:dimdrop} for an illustration of this idea.

It follows that one cannot expect (\ref{eq:sumcantorcompactdim})
to hold for {\em all} $a$, even when $K$ is also a central Cantor
set. In \cite{peresschlag} (as part of a more general framework)
the dimension of the set of exceptions in
(\ref{eq:sumcantorcompactdim}) was estimated, but whether the set
of exceptions (for $K=C_b$) is countable, was unknown until now.
Ero\v{g}lu \cite{eroglu-sums} proved that if $\dim(C_a)+\dim(C_b)
\le 1$, then the Hausdorff measure of the sum set $C_a+C_b$ in the
dimension $\dim(C_a)+\dim(C_b)$ is zero.

In a different direction, Moreira proved a result on the dimension
of the sum of dynamically-defined, non-linear Cantor sets. This
result was stated in \cite{moreira-hungarica}, but the proof
sketched there is incorrect; see \cite{shmerkin-moreira} for a
correct proof based on Moreira's ideas.

%\begin{theorem} (Moreira)
%Let $K$ (resp. $K'$) be the attractor of an i.f.s. $\{ f_1,\ldots,
%f_n\}$ (resp. $\{ g_1,\ldots, g_{n'}\}$), where $f_i, g_i$ are
%$C^{1+\e}$ diffeomorphisms of the unit interval $I$ into itself.
%Assume that the following conditions hold:
%\begin{enumerate}
%\item (nonlinearity) $\left(f_i^{-1} f_j\right)'(x_0)\neq 0$ for some $0< i< j\le n$, $x_0\in
%I$.
%\item (incommensurability)
%\[
%\log(g_1(y_1))/\log(f_1(x_1)) \notin \mathbb{Q},
%\]
%where $x_1$ (resp. $y_1$) is the fixed point of $f_1$ (resp.
%$g_1$).
%\end{enumerate}
%Then
%\[
%\dim(K+K') = \min(\dim(K)+\dim(K'),1).
%\]
%\end{theorem}

We remark that Moreira's Theorem has explicit checkable
conditions, but it applies to nonlinear constructions; on the
other hand the result of Peres and Solomyak applies to (affine)
self-similar sets but it is an almost-everywhere type of result
and gives no information about specific cases.

Theorem \ref{th:sumcantor} fills this gap by determining the exact
set of exceptions to (\ref{eq:sumcantorcompactdim}) in the case
where $K$ is also a central Cantor set. In fact, our result holds
for sums of general self-similar sets; see Theorem
\ref{th:sumselfsimilar} below.

Our research was also inspired by some conjectures of H.
Furstenberg [personal communication], who in the eighties asked
about the validity of Theorem \ref{th:sumcantor} in the particular
case where $a^{-1}$ and $b^{-1}$ are integers which are not the
power of a common integer. Other conjectures of Furstenberg, which
are similar in spirit, remain open. For example, let $S$ be the
one-dimensional Sierpinski gasket; i.e.
\[
S = \left\{ \sum_{i=1}^\infty z_i 3^{-i}: z_i\in
\{(0,0),(1,0),(0,1) \}\right\}.
\]
Furstenberg conjectured that $P_\theta(S)$ has dimension $1$ but
measure zero for all $\theta$ such that $\tan(\theta)$ is
irrational. The part of the conjecture concerning measure was
proved by Kenyon \cite{kenyon-projecting} and generalized by
Lagarias and Wang in \cite{lagarias-wang}, but the dimension part
is still open. Some related conjectures and results can be found
in \cite{furstenberg}.

 Palis~\cite{palis-conjecture} conjectured that, generically, if
$K+K'$ has positive Lebesgue measure, then it has nonempty
interior. This conjecture motivated much of the research on this
topic, and it was eventually answered positively in
\cite{moreirayoccoz}. Although our work is not directly related to
Palis' conjecture (indeed, we focus on the case where
$\dim(K)+\dim(K')<1$, and in particular $K+K'$ has zero measure),
it can be seen as a continuation of the same line of research. We
remark that Palis' conjecture restricted to self-similar sets is
still open; i.e. it is not known whether sums of self-similar sets
generically are either of zero measure or have nonempty interior.

The topological structure of $C_a+C_b$ when
$\dim(C_a)+\dim(C_b)>1$ was investigated in \cite{MendesOliveira}
and \cite{CabrelliHareMolter}; the condition $\log b/\log
a\notin\mathbb{Q}$ also arises naturally in this context.

\subsection{Further results}

The irrationality condition in Theorem \ref{th:sumselfsimilar} is
equivalent to $\{ \log|r_i|\} \cup \{\log|r'_i|\}$ not being an
arithmetic set (a subset of $\mathbb{R}$ is arithmetic if it is
contained in some lattice $\alpha\mathbb{N}$). Our next result
shows that, under the assumption that the Hausdorff dimension
equals the so-called similarity dimension, Theorem
\ref{th:sumselfsimilar} is sharp. We remark that this assumption
is weaker than the well-known open set condition, see e.g.
\cite{falconer2}.

\begin{theorem} \label{th:exceptions}
Let $K$ (resp. $K'$) be the attractor of the i.f.s. $\{ r_i x +
t_i \}_{i=1}^n$(resp. $\{r'_i x + t'_i\}_{i=1}^{n'}$), where $0<
r_i, r'_i<1$, and assume that the Hausdorff dimension of $K$
(resp. $K'$) is given by the only $\beta>0$ (resp. $\beta'>0$)
verifying $\sum_{i=1}^n r_i^\beta = 1$ (resp. $\sum_{i=1}^{n'}
(r'_i)^{\beta'}=1$). Suppose that the set $\{ \log |r_i| \} \cup
\{ \log |r'_i|\}$ is arithmetic. Then there exists some $s>0$ such
that
\begin{equation} \label{eq:dimdrop}
\dim(K+s K') < \dim(K) + \dim(K').
\end{equation}
\end{theorem}

Using self-similarity, it is easy to show that if there is one $s$
such that (\ref{eq:dimdrop}) holds, then there are infinitely
many. On the other hand, according to Marstrand's projection
theorem, the set of parameters $s$ such that a dimension drop
(\ref{eq:dimdrop}) occurs, has zero length. Furstenberg
conjectured that this set is in fact countable, at least in the
case where $K$ and $K'$ are central Cantor sets, but our methods
do not seem to yield progress on this problem.

A modification of the proof of Theorems \ref{th:sumcantor} and
\ref{th:sumselfsimilar} yields a new result on the projections of planar
self-similar sets generated by an iterated function system that includes a scaled
irrational rotation:

\begin{theorem} \label{th:irrationalrotation}
Let $\{ A_i x + d_i \}_{i=1}^n$ be a self-similar i.f.s. in
$\mathbb{R}^2$ with attractor $E$. Suppose that each linear map
$A_i$ is written as $\zeta_i R_{\theta_i} O_i$, where
$|\zeta_i|<1$, the map $R_{\theta_i}$ is the rotation by angle
$\theta_i$ and $O_i$ is either the identity or a reflection about
the $x$-axis. Suppose that the group $G$ generated by $\{
R_{\theta_i} O_i\}_{i=1}^n$ is such that the set $\{ \theta:
R_\theta\in G\}$ is dense in $[0,\pi)$. Then
\[
\dim(P_\xi(E)) = \min(\dim(E),1) \quad \textrm{ for all }
\xi\in [0,\pi).
\]
\end{theorem}

\begin{figure} \label{fig:ssrotations}
\begin{center}
\includegraphics[width=0.75\textwidth]{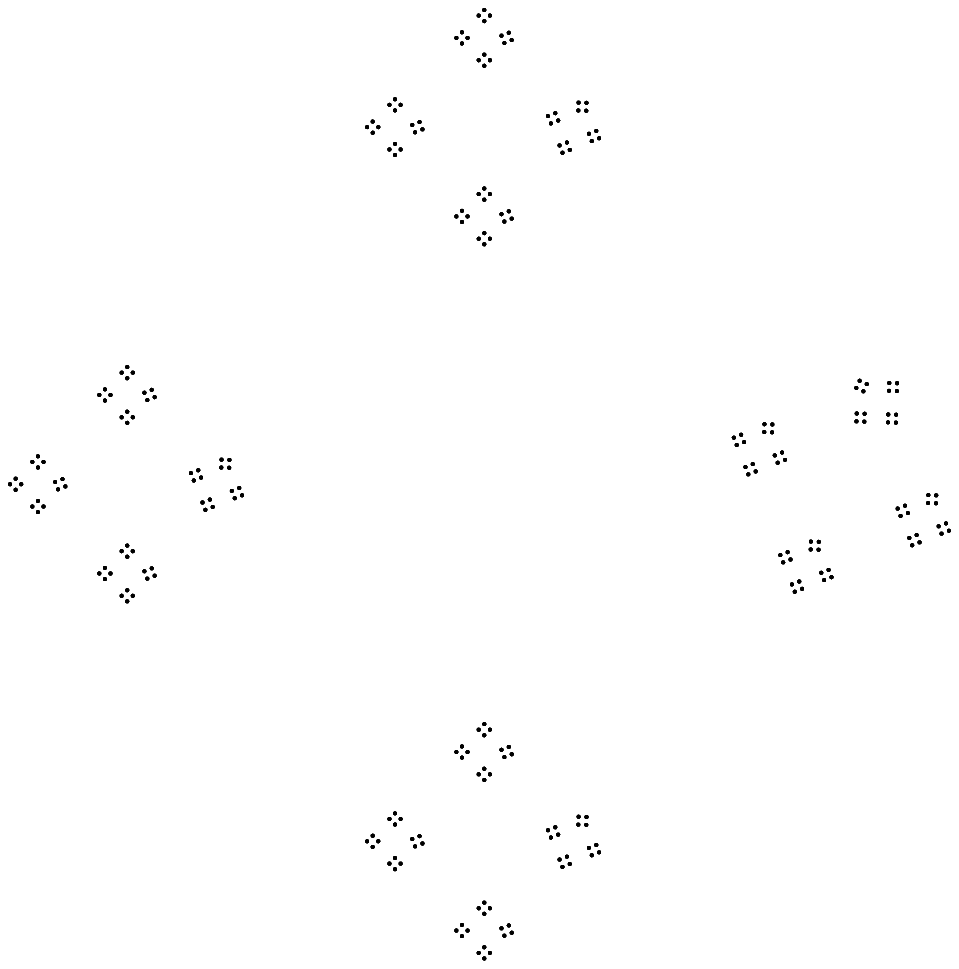}
\end{center}
\textbf{Figure 2}. The attractor of an i.f.s. to which Theorem
\ref{th:irrationalrotation} applies. The dimension of this
self-similar set is $1$. Theorem \ref{th:irrationalrotation} says
that all orthogonal projections are $1$-dimensional as well.
\end{figure}

Figure \ref{fig:ssrotations} depicts one of the self-similar sets
to which Theorem \ref{th:irrationalrotation} applies. We remark
that Ero\v{g}lu \cite{eroglu-projections} proved that, under the
assumptions of Theorem~\ref{th:irrationalrotation}, if the i.f.s.\ satisfies the open set condition and
$\gamma=\dim(E)\le 1$, then all projections $P_\xi(E)$ have
zero $\gamma$-dimensional Hausdorff measure. Thus, in this case, the
projections are smaller than the original set $E$ in measure, but
not in dimension.

As part of the proofs of Theorems \ref{th:sumselfsimilar} and
\ref{th:irrationalrotation}, we will need the following
proposition, which as far as we know is not in the literature.

\begin{prop} \label{prop:approxhomogeneous}
Let $\{ A_i x + d_i \}_{i=1}^n$ be a self-similar i.f.s. on
$\mathbb{R}^p$, where $p$ is either $1$ or $2$, satisfying the
open set condition, and let $E$ be its attractor.

Then for all $\e>0$, there exists an i.f.s. of the form $\{ A x +
z_i \}_{i=1}^N$ with attractor $H$, such that $H\subset E$, and
\[
\dim(H) \ge \dim(E) -\e.
\]
\end{prop}

\section{A discrete Marstrand projection theorem} \label{sec:preliminaries}

We prove a discrete version of Marstrand's Theorem on projections.
Similar results were obtained by Moreira (personal communication,
see \cite{shmerkin-moreira}) and Rams, in the context of
intersection numbers, see \cite{rams-pacific} and references
therein.

\begin{prop} \label{prop:discretemarstrand}
Given constants $A>1, A_1, A_2>0$ and $0<\gamma<1$, there exists a
constant $\delta$ such that the following holds:

Fix $\rho>0$. Let $\Q$ be a collection of disjoint closed convex
subsets of the unit disk such that each element contains a disk of
radius $A^{-1}\rho$ and is contained in a disk of radius $A\rho$.
Suppose that $\Q$ has cardinality at least $\rho^{-\gamma}/A_1$,
yet any disk of radius $\ell\in (\rho,1)$ intersects at most
$A_2(\ell/\rho)^\gamma$ elements of $\Q$. Then for any $\e>0$
there exists a set $J\subset [0,\pi]$ with the following
properties:
\begin{enumerate}
\item $\leb([0,\pi]\backslash J) \le \e$.
\item If $\theta\in J$, then there exists a subcollection
$\Q_1$ of $\Q$ of cardinality at least $\e \delta|\Q|$ such that
the orthogonal projections of the sets in $\Q_1$ onto a line with
direction $\theta$ are all disjoint and $\rho$-separated.
\item $J$ is a finite union of open intervals.
\end{enumerate}
\end{prop}
\textit{Proof}. In the course of the proof, $c$ denotes a
universal constant, and $A_i$ denote constants which depend only
on $A, A_1, A_2$ and $\gamma$. Let $E$ be the union of all
elements of $\Q$, and let $\mu$ be normalized Lebesgue measure on
$E$. Consider the Riesz energy
\[
I_1(\mu)=\int\int|z-w|^{-1} d\mu(z)d\mu(w).
\]
We claim that
\begin{equation} \label{eq:energybound}
I_1(\mu) \le A_4 \rho^{\gamma-1}.
\end{equation}
This is standard; we provide a proof below for completeness.

Assuming (\ref{eq:energybound}), we apply Theorem 9.9 in
\cite{mattila} with $n=2, m=1$, to obtain
\[
\int_0^\pi |P_\theta E|^{-1} d\theta \le cI_1(\mu) \le  A_5
\rho^{\gamma-1},
\]
where $P_\theta E$ are the projections of $E=\supp(\mu)$ onto a
line making angle $\theta$ with the $x$-axis, and $|P_\theta E|$
their length. Therefore, letting
\[
J_* = \{ \theta: |P_\theta E| > \e A_5^{-1} \rho^{1-\gamma} \}
\]
we have $\leb([0,\pi]\backslash J_*) \le \e$.

Let $J=J(\delta)$ be the set of all $\theta \in [0,\pi)$ with the following property:  there exists a subcollection
$\Q_1= \Q_1(\theta) $ of $\Q$ of cardinality at least $\e \delta|\Q|$, such that
the orthogonal projections of the sets in $\Q_1$ onto a line with
direction $\theta$ are all disjoint and $\rho$-separated.

Due to the convexity of the elements of $\Q$, the set
$J=J(\delta)$ is a finite union of open intervals. We claim  that
$J_* \subset J(\delta)$ for $\delta=(5A+5)^{-1}A_5^{-1}$. To prove
this, consider $\theta\in J_*$. By the definition of $J_*$,  we
can find in $P_\theta E$ at least $\delta \e \rho^{-\gamma}$
points $x_j$ that are $3(A+1)\rho$-separated. Choose for each of
these points $x_j$ one element in $\Q$ with projection that
contains $x_j$. This yields a family $\Q_1 \subset \Q$, of
cardinality at least $\delta \e \rho^{-\gamma}$, such that the
projections $\{ P_\theta(Q) : Q \in \Q_1\}$  are at least
$(A+1)\rho$-separated, as desired.

\medskip

\smallskip

\noindent\textit{Proof of (\ref{eq:energybound})}. We will need the
following properties of $\mu$:
\begin{enumerate}
\item[(i)] $A_2^{-1} A^{-4}\rho^\gamma \le \mu(Q) \le A_1 A^4 \rho^\gamma$ for any element $Q$
of $\Q$.
\item[(ii)] $\mu(D) \le A_7 \ell^\gamma$ for any disk $D$ of
radius $\ell\in [\rho,1]$.
\item[(iii)] $\mu \le A_8 \rho^{\gamma-2}\leb_2$, where $\leb_2$
denotes two-dimensional Lebesgue measure.
\end{enumerate}
To begin, observe that if $Q_1, Q_2$ are elements of $\Q$, then
$\mu(Q_1)/\mu(Q_2)\le A^4$. Since $\Q$ has between $
\rho^{-\gamma}/A_1$ and $A_2 \rho^{-\gamma}$ elements, (i)
follows. If $D$ is a disk of radius $\ell$, where $\rho\le\ell\le
1$, then $D$ intersects at most $A_2(\ell/\rho)^\gamma$ elements
of $\Q$, whence (ii) follows from (i) with $A_7 = A_1 A^4 A_2$.
Next, observe that
\[
\leb_2(E) \ge (A_1^{-1}\rho^{-\gamma})(A^{-2}\rho^2) = A_1^{-1}
A^{-2} \rho^{2-\gamma},
\]
and from this we deduce that, for any Borel set
$B\subset\mathbb{R}^2$,
\[
\mu(B) = \mu(B\cap E) = \leb_2(E)^{-1}\leb_2(B\cap E) \le A_1 A^2
\rho^{\gamma-2}\leb_2(B),
\]
which shows that (iii) holds with $A_8 = A_1 A^2$.

Now fix $w\in E$. Using (iii) we estimate
\begin{equation} \label{eq:energyest1}
\int_{|z-w|\le \rho} |z-w|^{-1} d\mu(z) \le A_8 \rho^{\gamma-2}
\int_{|u|\le\rho} |u|^{-1} du = 2\pi A_8 \rho^{\gamma-1}.
\end{equation}
On the other hand, for $j\ge 1$ we can apply (ii) to obtain
\[
\mu(\{z:|z-w|\le e^j\rho\}) \le A_7 \left(e^j\rho\right)^\gamma,
\]
and from here we deduce
\begin{equation} \label{eq:energyest2}
\int_{e^{j-1}\rho < |z-w|\le e^j\rho} |z-w|^{-1}d\mu(z) \le  A_7
\left(e^j\rho\right)^\gamma e^{1-j}\rho^{-1} = e A_7
(e^j\rho)^{\gamma-1}.
\end{equation}
Adding (\ref{eq:energyest1}) and the sum of (\ref{eq:energyest2})
over all $j\in\mathbb{N}$, we get
\[
\int |z-w|^{-1} d\mu(z)  \le A_4 \rho^{\gamma-1}.
\]
Finally, integrating over $w$ yields (\ref{eq:energybound}), and
the proof is complete. \qed

\section{Proof of Theorem \ref{th:sumcantor}}
\label{sec:proofmaintheorem}

We first discuss the main ideas of the proof of Theorem
\ref{th:sumcantor}; full details follow below.
Recall that Marstrand's theorem says that
\[
\dim(C_a+ e^t C_b) = \min(\dim(C_a)+\dim(C_b),1)
\]
for almost every $t\in\mathbb{R}$ (for reasons which will become
apparent later, it will be useful to work in a logarithmic scale).
Proposition \ref{prop:discretemarstrand} gives a discrete version
of this: we decompose $C_a\times C_b$ into cylinder rectangles of
size approximately $\rho\times\rho$. For $t\in\mathbb{R}$ consider
the projection mapping $\Pi_{e^t}$. It will follow from
Proposition \ref{prop:discretemarstrand} that there is a set $F_0$
of ``good'' values, whose complement has very small measure, such
that for $t\in F_0$ a large number of rectangles have disjoint
$\Pi_{e^t}$ projections.

Because of the homogeneity of $C_a, C_b$, all these cylinders are
translates of each other, and they are pairwise disjoint. The idea
now is to decompose these first cylinders into subcylinders of
size approximately $\rho^2\times\rho^2$ and consider the set of
``good'' scales $F_1$ associated to them; because the cylinders
are a (anisotropic) rescaling of $C_a\times C_b$, $F_1$ is of the
form $F_0 + \alpha_1$ for some $\alpha_1\in\mathbb{R}$. We
continue this process inductively, to obtain sets of good scales
$F_n$ after $n$ steps; we have $F_{n+1} = F_n + \alpha_n$ for some
sequence $\{\alpha_n\}$. Notice that there is some latitude in
choosing this sequence; the main idea is to exploit this freedom
to pick the $\alpha_n$ in such a way that $t\rightarrow
t+\alpha_n$ is essentially an irrational rotation on the circle.
By construction the sets $F_i$ are robust: they can be taken to be
a finite union of intervals. Then Weyl's equidistribution
principle implies that for any specific angle $\theta$, there are
many cylinders at many levels (more precisely, at a set of levels
of well-defined density close to $1$) whose images under
$\Pi_{e^t}$ are disjoint. From here a standard application of the
mass distribution principle gives the desired lower bound on the
dimension.

We remark that the idea of using the discrete version of
Marstrand's theorem and then iterating the process on each
cylinder is due to Moreira. However, while he uses a recurrence
result which requires a nonlinearity assumption, we go in a
different direction by using Weyl's equidistribution principle.
Furstenberg used Weyl's Theorem in a related setting in
\cite{furstenberg}.

We now proceed to the details of the proof. We will establish
a slightly more general version which will be needed to prove
Theorem \ref{th:sumselfsimilar}.

\begin{theorem} \label{th:sumhomogeneous}
Let $K$ be the attractor of the i.f.s. $\{ rx+t_i\}_{i=1}^n$, and
let $K'$ be the attractor of the i.f.s. $\{r'x+t'_i\}_{i=1}^{n'}$.

Let $I, I'$ be the convex hulls of $K, K'$ respectively. Assume
that the intervals $r I+t_i$ are pairwise disjoint, as are the
intervals $r' I' + t'_i$.

Assume further that $\log(|r|)/\log(|r'|)$ is irrational.

Then
\begin{equation} \label{eq:equalityforalldilations}
\dim(K+ K') = \min(\dim(K)+\dim(K'),1),
\end{equation}
\end{theorem}

The central Cantor sets $C_a$ satisfy the hypothesis of Theorem
\ref{th:sumhomogeneous}, whence Theorem \ref{th:sumcantor} is an
immediate consequence of it.

\textit{Proof of Theorem \ref{th:sumhomogeneous}}.

\textbf{Notation and remarks}. The upper estimate
\[
\dim(K+K') \le \min(\dim(K) + \dim(K'),1).
\]
is well known, so we only need to prove that the opposite
inequality also holds.

If $\dim(K)+\dim(K')\ge 1$, then for any $\e>0$ we can, by
iterating the original i.f.s. and throwing away some maps, find a
subset $\widetilde{K}$ of $K$ such that the hypothesis of the
Theorem still apply to the pair $(\widetilde{K}, K')$, and
\[
1-\e < \dim(\widetilde{K})+\dim(K') < 1.
\]
Thus we can assume without loss of generality that
\[
\gamma:=\dim(K)+\dim(K')<1.
\]

Dilating both $K, K'$ by the same factor and translating either of
$K$ or $K'$ does not affect either the hypothesis nor the result;
therefore we can also assume that $I=[0,1]$ and $I'=[0,e^\tau]$
for some $\tau\in\mathbb{R}$. Therefore it is enough to show that
\[
\dim(K+ e^\tau K') \ge \dim(K) + \dim(K') - \eta
\]
for all $\eta>0$, under the additional assumption that the unit
interval is the convex hull of both $K$ and $K'$. We will fix
$\tau$ for the rest of the proof.

Let $\Sigma^*$ be the family of all finite words
$u=(u_1,\ldots,u_k)$ with $u_i\in\{1,\ldots,n\}$; we define
$(\Sigma')^*$ analogously. As customary, $|u|$ will denote the
length of $u$. Concatenation of words will be denoted by
juxtaposition; $1^j$ denotes the word consisting of $j$
consecutive ones.

Let $f_i(x) = rx+t_i$, $f'_i(x)=r' x+t'_i$. For $u\in \Sigma^*$
,$u'\in(\Sigma')^*$ we will write
\[
f_u = f_{u_1}\cdots f_{u_k},\quad f'_{u'} = f'_{u'_1}\cdots
f_{u'_{k'}},
\]
where $|u|=k, |u'|=k'$. Also let $I(u)=f_u(I)$ and
$I'(u')=f'_{u'}(I)$. Finally, write $Q(u,u') = I(u)\times I'(u')$.
We underline that $f'$ here is not a derivative.

Recall that $P_\theta$ is the orthogonal projection onto a line
making angle $\theta$ with the $x$-axis. The projection mappings
$(x,y)\rightarrow x+ty$ will be denoted by $\Pi_t$. A trivial but
useful observation is that $P_\theta(\Lambda)$ and
$\Pi_{\tan\theta}(\Lambda)$ are the same up to affine equivalence.

\bigskip
\noindent

\noindent\textbf{Application of the discrete version of Marstrand's
Theorem}.

\medskip

In the course of the proof, $c$ will denote a constant which
depends only on $K$ and $K'$; its value may be different in each
line.

Given $k\ge 0$, consider the families
\begin{eqnarray*}
 \Q_k &  = & \{ Q(u,u'): |u|=k, |u'|=k' \},\\
\widetilde{\Q}_k & = & \{ Q(u,u'): |u|=k, |u'|=k'+1 \},
\end{eqnarray*}
where $k'$ is the largest integer such that $r^k < (r')^{k'}$.
Note that rectangles in $\Q_k$ have size $r^k \times M_k r^k$,
where $1<M_k<1/r'$, and rectangles in $\widetilde{\Q}_k$ have size
$r^k \times M_k r' r^k$.

Observe that $|\Q_j| > c r^{-j \gamma}$. This is well known; see
for example \cite[Proposition 7.4]{falconer3} and the remarks
preceding it. Moreover, if $j<k$, then any $Q\in \Q_j$ contains at
most $c r^{(j-k)\gamma}$ elements of $\Q_k$. Therefore if $S$ is
any square of side $r^j$, then $S$ intersects a uniformly bounded
number of elements of $\Q_j$, and combined with the previous
remark we get that
\[
| \{Q\in \Q_k : Q\subset S\} | \le c r^{(j-k)\gamma}.
\]
It follows from these observations that $\Q_k$ verifies the
conditions of Proposition \ref{prop:discretemarstrand} for some
$A, A_1, A_2$ which depend only on $K,K'$, and in particular are
independent of $k$.

Now pick some large integer $m$ and some small $\e>0$, and apply
Proposition \ref{prop:discretemarstrand} to $\Q_m$ to obtain a set
$J$ of ``good'' angles. Note that because every rectangle of
$\widetilde{\Q}_m$ is contained in a rectangle of $\Q_m$, property
(2) of $J$ in Proposition \ref{prop:discretemarstrand} is valid
for $\widetilde{\Q}_m$ as well.

Write $\alpha=\log(M_m)$, $\beta = \log(1/r')$. A crucial fact is
that $\alpha/\beta$ is irrational; this follows easily from the
irrationality assumption. In particular, the map
\begin{equation} \label{eq:defrotation}
\rot(x) = x+\alpha\,(\textrm{mod} \beta)
\end{equation}
is uniquely ergodic on $[0,\beta)$ endowed with normalized
Lebesgue measure. We observe that $\rot^k(0) = \log(M_{mk})$; this
fact will be useful later.

We will show that
\[
\dim\left(K + e^{\tau}K'\right) \ge \gamma-\eta,
\]
where $\eta$ depends on $m$ and $\e$, and can be made arbitrarily
small by letting $m\rightarrow\infty, \e\rightarrow 0$.

If $\theta\in [0,\pi)$, then the orthogonal projection $P_\theta$
can be identified with the map $\Pi_s(x,y) = x+sy$, where
$s=\tan(\theta)$. Since the map $\theta \rightarrow
\log(\tan(\theta))$ is a diffeomorphism on any compact subset of
the set of directions with positive and finite slope, it follows
from Proposition \ref{prop:discretemarstrand} that there exists
$\widetilde{F}\subset [\tau,\tau+\beta)$ and constants $L,
\delta_1$ (independent of $m$ and $\e$) such that
\begin{enumerate}
\item[(i)] $\leb([\tau,\tau+\beta)\backslash \widetilde{F}) \le L\e$.
\item[(ii)] If $t\in \widetilde{F}$, then there exist $\D=\D(t)\subset \Q_m ,\widetilde{\D}=\widetilde{\D}(t)\subset  \widetilde{\Q}_m$ such that $|\D|> \delta_1
\e r^{-m\gamma}$ and the family
\[
\left\{ \Pi_{e^t}(Q) : Q\in \D \right\}
\]
is $r^m$-separated; and analogous assertions hold for
$\widetilde{D}$. '
\item[(iii)] $\widetilde{F}$ is a union of finitely many open intervals.
\end{enumerate}
Finally, let $F = \widetilde{F}-\tau$.

\bigskip
\noindent

\textbf{Inductive construction}. We will now use properties (i)
and (ii) to inductively construct a tree $\R$ with vertices
labeled by cylinder rectangles. The set of rectangles of level $j$
will be denoted by $\R_j$. We will prove that the tree has the
following properties:
\begin{enumerate}
\item[(A)] If $Q$ is the parent of $Q'$, then $Q'\subset Q$.
\item[(B)] If $Q\in \R_j$, then $Q$ has size $r^{mj}\times
\exp(\rot^j(0)) r^{mj}$.
\item[(C)] The family $\{ \Pi_{e^\tau}(Q) : Q\in \R_j\}$ is disjoint and
$r^{jm}$-separated.
\item[(D)] All rectangles of level $j$ have the same number of children $C_j$. Moreover,
\begin{eqnarray*}
\rot^j(0) \notin F & \Longrightarrow & C_j =1,\\
\rot^j(0) \in F & \Longrightarrow & C_j \ge \delta_1\e
r^{-m\gamma}.
\end{eqnarray*}
\item[(E)] Each $Q\in\R_j$ is of the form $Q(v,v')$ for some
$v\in\Sigma^*, v'\in(\Sigma')^*$.
\end{enumerate}

Let the root of the tree be the unit square. Define $\R_1$ as
follows: if $\tau\in \widetilde{F}$, then $\R_1=\D(t)$; otherwise,
$\R_1 = \{ \widehat{Q}\}$, where $\widehat{Q}$ is any element of
$\Q_m$ (say, $\widehat{Q}=Q(1^m,1^{m'}$)). Property (B) follows
from the definition of $\alpha$. Properties (A), (C), (D) and (E)
are also clear.

Now assume that $\R_j$ has been defined verifying (A)-(E). We will
consider two cases:
\begin{itemize}
\item[(I)] $\rot^j(0) + \alpha < \beta$.
\item[(II)] $\rot^j(0) + \alpha > \beta$.
\end{itemize}
Let $\sigma=0$ if (I) holds, and $\sigma=1$ otherwise. Observe
that
\begin{equation} \label{eq:nextrotation}
\rot^{j+1}(0) = \rot^j(0)+\alpha-\sigma\beta.
\end{equation}
If $P, Q$ are rectangles define their product $Q\cdot P$ as $f_Q
f_P(Q_0)$, where $Q_0$ is the unit square and $f_P, f_Q$ are the
orientation preserving linear maps mapping $Q_0$ onto $P, Q$
respectively. Now given $Q\in \R_j$ we will define its set of
offspring $C(Q)$ as follows: if $\rot^j(0)\in F$, let
\[
C(Q) = \left\{ \begin{array}{cc}
       \{  Q\cdot P : {P \in \D(\rot^j(0)+\tau)} \}  & \textrm{ if case (I) holds } \\
                 \{  Q\cdot P : {P \in \widetilde{\D}(\rot^j(0)+\tau)} \} & \textrm{ if case (II) holds }  \\
               \end{array}\right..
\]
Otherwise, if $\rot^j(0)\notin F$, then let $C(Q) = Q\cdot Q(1^m,
1^{m'+\sigma})$.

Property (A) is clear. Property (B) follows from
(\ref{eq:nextrotation}) since each element of $\R_{j+1}$ has size
\[
r^{mj} r^m \times r^{mj} \exp(\rot^j(0)) M_m (r')^{\sigma}r^m =
r^{m(j+1)} \times \exp(\rot^j(0)+\alpha-\sigma\beta) r^{m(j+1)} .
\]
If $\rot^j(0)\notin F$, then property (C) is trivial, so we will
assume that $\rot^j(0)\in F$ or, in other words,
$\rot^j(0)+\tau\in \widetilde{F}$. It is enough to show that if
$S_1, S_2$ are children of the same $Q\in\R_j$, then their
projections are $r^{m(j+1)}$ separated. Moreover, by translating
$Q$ if needed we can assume that $Q=r^{mj}(I \times \exp(\rot^j
(0))I)$ where $I$ is the unit interval. Let $S_i = Q \cdot P_i$,
where $P_i = J_i \times J'_i$. We have:
\begin{eqnarray*}
\Pi_{e^{\tau}}(S_i) &  = & \Pi_{e^{\tau}}(Q\cdot P_i)\\
& = & r^{jm} \Pi_{e^{\tau}}( (I \times \exp(\rot^j(0)) I) \cdot
(J_i \times J_i'))\\
& = & r^{jm} (J_i \times \exp(\rot^j(0)+\tau) J_i')\\
& = & r^{jm} \Pi_{\exp(\rot^j(0)+\tau)} P_i.
\end{eqnarray*}
By property (ii) of $\widetilde{F}$ the family
$\left\{\Pi_{\exp(\rot^j(0)+\tau)} P_i\right\}$ is
$r^m$-separated; it follows that the family
$\left\{\Pi_{e^{\tau}}(S_i)\right\}$ is $r^{(j+1)m}$-separated, as
desired.

Property (D) follows immediately from Property (ii) above.
Finally, Property (E) is a direct consequence of the construction.

\bigskip
\noindent

\textbf{Application of Weyl's equidistribution principle}.

We recall that $\widetilde{F}$, and hence $F$, is a finite union
of open intervals. Therefore we can apply Weyl's equidistribution
principle to get
\begin{equation} \label{eq:weyl}
\lim_{j\rightarrow\infty} \frac{1}{j}\left|\left\{ i<j : \rot^i(0)
\in F\right\}\right| = \frac{\leb(F\cap
[0,\beta))}{\leb([0,\beta))} \ge 1 - (L/\beta) \e.
\end{equation}

Let
\[
E_\tau = \bigcap_{j=1}^\infty \bigcup_{Q\in \R_j} \Pi_{e^\tau}(Q).
\]
By property (A) this set is a countable intersection of nested
nonempty compact sets, so $E_\tau$ is compact and nonempty. It
follows from property (E) that $E_\tau\subset K+e^\tau K'$.

We will estimate the dimension of $E_\tau$ in a standard way, by
means of a Frostman measure. Let $\mu_\tau$ be the probability
measure which assigns the same mass $|\R_j|^{-1}$ to all intervals
$\Pi_{e^\tau}(Q)$ for $Q\in\R_j$. This measure is well defined
thanks to properties (C) and (D), and is supported on $E_\tau$ by
definition.

Let $x\in E_\tau$. Because of (C) the interval
$(x-r^{jm},x+r^{jm})$ intersects exactly one interval
$\Pi_{e^\tau}(Q)$ with $Q\in\R_j$, so we have
\[
\mu_\tau(x-r^{jm},x+r^{jm}) \le \mu_\tau(\Pi_{e^\tau}(Q)) =
|\R_j|^{-1}.
\]
On the other hand, using (D) once again we get that
\[
\log|\R_j| \ge \left|\left\{ i<j : \rot^i(0) \in F\right\}\right|
\log\left(\delta_1\e r^{-m\gamma}\right).
\]
Hence we deduce from (\ref{eq:weyl}) that if $j$ is large enough,
then
\[
\log\mu_\tau(x-r^{jm},x+r^{jm}) \le j(1-(2L/\beta)\e)
\log\left(\delta_1^{-1}\e^{-1} r^{m\gamma}\right).
\]
Thus from the mass distribution principle (see e.g.
\cite[Proposition 2.3]{falconer3}) we conclude that
\[
\dim(K+e^\tau K')\ge\dim(E_\tau) \ge  \frac{(1-(2L/\beta)\e)
\log\left(\delta_1^{-1}\e^{-1} r^{m\gamma}\right)}{m\log r}.
\]

The right-hand side can be made arbitrarily close to $\gamma$ by
letting $m\rightarrow\infty$, $\e\rightarrow 0$. Therefore the
Hausdorff dimension of $K + e^\tau K'$ must be at least $\gamma$,
and this completes the proof. \qed

\section{Remaining proofs} \label{sec:remainingproofs}

\medskip

\textit{Proof of Theorem \ref{th:sumselfsimilar}}. Let $K$ and
$K'$ be as in the statement of the theorem. We can assume without
loss of generality that
$\log(|r_1|)/\log(|r'_1|)\notin\mathbb{Q}$.

Fix $\e>0$. Pick $\delta$ small enough so that for any covering
$\{ I_j\}$ of $K$ by intervals of diameter at most $3\delta$, we
have
\[
\sum_j |I_j|^{\dim(K)-\e} > 1.
\]

Let $I$ be the convex hull of $K$. Consider a covering of $K$ by
cylinder intervals $f_u(I)$ of diameter between $r_* \delta$ and
$\delta$, where $r_* = \min_i |r_i|$. Pick a maximal pairwise
disjoint collection $\{ f_{u_j}(I) \}_{j=1}^N$ among these
cylinders. By maximality, if $I_j$ is the interval with the same
center as $f_{u_j}(I)$ and length $|f_{u_j}(I)|+2\delta$, then $\{
I_j\}_{j=1}^N$ is a covering of $K$. Then we have
\[
\sum_{j=1}^N (3\delta)^{\dim(K)-\e} \ge \sum_{j=1}^N
|I_j|^{\dim(K)-\e}
> 1.
\]
It follows that $N\ge c \delta^{\e-\dim(K)}$, where
$c=3^{\e-\dim(K)}$. Let $\widetilde{K}$ be the attractor of the
i.f.s. $\{f_{u_j}\}_{j=1}^N$. We have that $\widetilde{K}\subset
K$ and
\[
\dim(\widetilde{K}) \ge \frac{\log N}{|\log(r_*\delta)|} \ge
\frac{\log c +
(\dim(K)-\e)|\log\delta|}{|\log(r_*)|+|\log\delta|}.
\]
Hence by taking $\delta$ very small we can ensure that
$\dim(\widetilde{K}) > \dim(K)-2\e$. Moreover, all the intervals
$f_{u_j}(I)$ are disjoint. In the same way we obtain an
appropriate set $\widetilde{K}'\subset K'$ . Therefore we can
assume, without loss of generality, that the families
\[
\{ r_i I + t_i \}_{i=1}^n,\quad \{ r'_i I' + t'_i \}_{i=1}^{n'}
\]
are disjoint for some intervals $I$, $I'$ (one can take these
intervals to be the convex hulls of $K$ and $K'$). At this point
the irrationality condition may fail to hold.

We now apply Proposition \ref{prop:approxhomogeneous} to $K$ and
$K'$, to obtain iterated function systems $\{ g_i(x) \}_{i=1}^N$,
$\{ g'_i(x) \}_{i=1}^{N'}$ with attractors $H, H'$, where
\[
g_i(x) = \rho x + z_i, \quad g'_i(x) = \rho' x + z'_i,
\]
the dimensions of $H, H'$ are very close to those of $K, K'$, and
$H\subset K, H'\subset K'$. If $\log|\rho|/\log|\rho'|$ is
irrational we can apply Theorem \ref{th:sumhomogeneous} to $H+H'$,
and we are done. If not, then consider the new i.f.s.
\[
\{ f_1 \circ g_i(x) \}_{i=1}^N,
\]
where $f_1(x) = r_1 x + t_1$. The dimension of the attractor
$\widetilde{H}\subset K$ is given by
\[
\frac{\log N}{-(\log|\rho|+\log|r_1|)},
\]
which, by picking $\rho$ small enough, can be made arbitrarily
close to the dimension of $H$ (and hence to the dimension of $K$).
Define $\widetilde{H}'$ analogously. Since we are assuming that
$\log(|\rho|)/\log(|\rho'|)\in\mathbb{Q}$, and that
$\log(|r_1|)/\log(|r'_1|)\notin \mathbb{Q}$, it follows that
$\widetilde{H}, \widetilde{H}'$ verify the irrationality
hypothesis and we can apply Theorem \ref{th:sumhomogeneous} to
$\widetilde{H}+\widetilde{H}'$. This completes the proof. \qed

\bigskip

\textit{Proof of Theorem \ref{th:exceptions}}.

By assumption there exist $\xi>0$ and integers $\{ a_i\}$, $\{
b_i\}$ such that $r_i=\xi^{a_i}$ and $r'_i=\xi^{b_i}$.

We first consider the case in which all the $a_i$ are equal to
some $a$ and all the $b_i$ are equal to some $b$. In this case
$K+K'$ is also self-similar and the proof is much easier; notice
that algebraically resonant central Cantor sets fall into this
case. By iterating the first i.f.s. $b$ times and the second
i.f.s. $a$ times, and replacing $\xi$ by $\xi^{ab}$, we can assume
that $a=b=1$. In this case the attractors are up to affine
equivalence explicitly given by
\[
K = \left\{ \sum_{i=1}^\infty d_i \xi^i : d_i \in \{t_1,\ldots,
t_n\}\right\},\quad K' = \left\{ \sum_{i=1}^\infty d_i \xi^i : d_i
\in \{t'_1,\ldots, t'_{n'}\}\right\}.
\]
Moreover, $\dim(K)=\log(n)/\log(1/\xi)$, and
$\dim(K')=\log(n')/\log(1/\xi)$. Therefore, letting $D=\{
t_1,\ldots, t_n\}$ and $D'=\{t'_1,\ldots, t'_{n'}\}$, we have
\[
K + s K' = \left\{ \sum_{i=1}^\infty d_i \xi^i : d_i\in D+ s D'
\right\},
\]
whence
\[
\dim(K+ s K') \le \frac{\log|D+s D'|}{\log(1/\xi)}.
\]
One can clearly take $s$ so that $|D+sD'|<|D||D'|$; for example,
let $s=(t_n-t_1)/(t'_{n'}-t'_1)$. For any such $s$, a dimension
drop occurs.

We now consider the general case. We will use the same notation as
in the proof of Theorem \ref{th:sumhomogeneous}. We claim that by
iterating the constructions and reordering we can assume that $a_1
= a_2 = b_1 = b_2$. Indeed, one can replace $f_i\, (i=1,2)$ by the
collection of maps $\{ f_i f_u : |u|=a_{3-i} b_1 b_2-1\}$, and
$f'_i\, (i=1,2)$ by  $\{ f'_i f'_u : |u|=b_{3-i} a_1 a_2-1\}$.
This operation does not change the attractors $K$ and $K'$. But
now $f_1 \circ f_1^{a_2 b_1 b_2-1}$ and $f_2 \circ f_2^{a_1 b_1
b_2-1}$ are maps in the first i.f.s. (they are obtained by taking
$i=1$ and $u$ a word consisting of all $1$s, and then $i=2$ and
$u$ a word consisting of all $2$s, respectively), and they both
have similarity ratio $\xi^{a_1 a_2 b_1 b_2}$. Likewise, there are
two maps in the second i.f.s. with this similarity ratio; this
proves the claim.

Denote the common value of $a_1, a_2, b_1, b_2$ by $\ell$. By a
further translation and dilation we can also assume that the
convex hull of both $K$ and $K'$ is the unit interval.

After these simplifications, notice that $Q(1,1)$ and $Q(2,2)$ are
squares of the same size, and we can choose $s$ so that
$\Pi_s(Q(1,1))=\Pi_s(Q(2,2))$. By our assumption that all the maps
are orientation-preserving and by self-similarity, it follows that
\[
\Pi_s(f_1(K)\times f'_1(K')) = \Pi_s(f_2(K)\times f'_2(K')).
\]
More generally, if for some $u, u'$ we have $r_u = r'_{u'}$, then
\begin{equation} \label{eq:coincidingprojections}
\Pi_s(f_{(u1)}(K)\times f'_{(u'1)}(K')) =  \Pi_s(f_{(u2)}(K)\times
f'_{(u'2)}(K')).
\end{equation}

Let $a$ be the g.c.d. of $\{a_i\}$, and $b$ the g.c.d. of
$\{b_i\}$. Write also $A=\max\{a_i\}, B=\max\{b_i\}$. By
\cite[Theorem 1.4.1]{kemeny-snell}, there is a smallest integer
$M_0$ with the following property: all multiples of $ab$ that are
greater than or equal to $M_0$ can be represented as a linear
combination of $\{a_i\}$ with positive integral coefficients, and
can also be  represented as a linear combination of $\{b_i\}$ with
positive integral coefficients. Pick some $M\ge M_0+\max(A,B)$
which is a multiple of $ab$.

Let $\Sigma = \{1,\ldots,n\}^{\mathbb{N}}$. Given
$\omega\in\Sigma$ let
\begin{eqnarray*}
S_j(\omega) & = & \sum_{i=1}^j a_{\omega_i};\\
n_k(\omega) & = & \min\{ j: S_j(\omega)\ge k M  \};
\\L_k(\omega) & = & S_{n_k}(\omega)- k M.
\end{eqnarray*}

All these definitions also apply to finite words of the
appropriate length. Recall that the Hausdorff dimension of $K$ is
given by the only $\beta>0$ such that
\[
\sum_{i=1}^n \xi^{\beta a_i} = 1.
\]
Endow $\Sigma$ with the Bernoulli measure $\mu$ for the
probability vector $\left(\xi^{\beta a_i}\right)_{i=1}^n$. Define
$\Sigma', S'_j, n'_k, L'_k, \beta', \mu'$ analogously.

Notice that $L_k$ can take the values $0, a,\ldots, (A/a-1)a$.
Because of the way $M$ was defined, there is a finite sequence $\{
j_m\}$ such that
\[
(k+1)M = S_{n_k}(\omega) + \sum a_{j_m}.
\]
Since $\mu$ is Bernoulli, this shows that there exist $\{p_i\}$,
independent of $k$, such that
\begin{equation} \label{eq:hittingprob}
\mu(S_{n_{k+1}} = (k+1)M |  L_k = i a   ) = p_i,\quad \textrm{ for
} i=0,1,\ldots A/a-1,
\end{equation}
and analogously for $\mu'$. Let $\nu=\mu\times\mu'$; we will treat
$\nu$ like a probability distribution according to which pairs of
sequences are drawn at random. Consider the events
\begin{eqnarray*}
\Omega_k &= & \left\{ S_{n_k}
=S'_{n'_k}=kM \right\}; \\
\Xi_{k,i} & = & \left\{ \omega(n_k+1)=i,
\omega'(n'_k+1)=i \right\} \quad (i=1,2);\\
\Theta_{k,i} & = & \Omega_k \cap \Xi_{k,i}.
\end{eqnarray*}
Letting $p=\min(\{p_i\},\{p'_i\})$, it follows from
(\ref{eq:hittingprob}) and independence that
\[
\nu(\Omega_k) \ge p^2, \quad \nu(\Xi_{k,i}) =
\xi^{\ell(\beta+\beta')},\quad \nu(\Theta_{k,i})\ge
\xi^{\ell(\beta+\beta')} p^2.
\]
Therefore
\begin{eqnarray*}
\nu(\Theta_{k+1,2}\cap\Theta_{k,2}^c) & > &
\nu(\Theta_{k+1,2}\cap \Theta_{k,1})\\
& = & \nu\left((\Theta_{k+1,2}|S_{n_k+1}=S'_{n'_{k'}+1}=k M+\ell ) \cap \Theta_{k,1}\right) \\
& \ge & \left(\xi^{\ell(\beta+\beta')} p^2\right)^2 =: q,
\end{eqnarray*}
where in the last inequality we used independence again. In
particular, we have showed that
$\nu(\Theta_{k+1,2}^c|\Theta_{k,2}^c)< 1-q$, and from this we
deduce that
\begin{equation} \label{eq:probessential}
\nu\left( \bigcap_{j=1}^k \Theta_{j,2}^c \right) < (1-q)^k.
\end{equation}
We now come back to the geometric picture. Consider the following
family:
\[
\mathcal{C}_k = \left\{ u\in\Sigma^*: S_{|u|-1}(u) < k M \le
S_{|u|}(u) \right\}.
\]
Notice that $\{ f_u(K): u\in\mathcal{C}_k\}$ is a covering of $K$;
define $\mathcal{C}'_k$ analogously. Let us call a pair
$(u,u')\in\mathcal{C}_k\times \mathcal{C}'_k$ {\em redundant} if
$(u,u')\in \Theta_{j,2}$ for some $j<k$; otherwise, let us call it
{\em essential}.

We claim that
\[
\mathcal{E}_k = \left\{ \Pi_s(Q(u,u')): (u,u')\in
\mathcal{C}_k\times \mathcal{C}'_k \textrm{ is essential }
\right\}
\]
is a covering of $\Pi_s(K\times K')$. Indeed, if $(u,u')$ is
redundant, then it can be decomposed as $(v_1 2 v_2, v'_1 2 v'_2)$
where $r_{v_1}=r'_{v'_1}=\xi^{jM}$. Hence it follows from
(\ref{eq:coincidingprojections}) that
\[
\Pi_s(K_u \times K'_{u'}) \subset \Pi_s(K_{(v_1 1)} \times
K'_{(v'_1 1) }) ,
\]
and if we delete $(u,u')$ from $\mathcal{C}_k\times
\mathcal{C}'_k$ we still get a covering of the projection.

Notice that each interval in $\mathcal{E}_k$ has length bounded by
$(s+1)\xi^{k M}$. Any essential pair $(u,u')$ has
$\nu$-probability at least
$\xi^{A\beta+B\beta'}\xi^{k(\beta+\beta')M}$ of occurring. On the
other hand, the probability of a pair in
$\mathcal{C}_k\times\mathcal{C}'_k$ being essential is, by
(\ref{eq:probessential}), no more than $(1-q)^{k-1}$. It follows
that there are at most
\[
\frac{1}{(1-q)\xi^{A\beta+B\beta'}}\left(\left(\xi^{-(\beta+\beta')M}\right)(1-q)\right)^k
\]
essential pairs. As $k\rightarrow\infty$ we get coverings of
$K+e^s K'$ by intervals of arbitrarily small diameter. In the
limit we can ignore the factors which do not depend on $k$, and we
conclude that
\[
\dim(K+ e^s K') \le
\frac{(\beta+\beta')M\log(1/\xi)+\log(1-q)}{M\log(1/\xi)} <
\beta+\beta',
\]
as desired. \qed

\bigskip

\textit{Proof of Theorem \ref{th:irrationalrotation}}. By
proceeding like in the proof of Theorem \ref{th:sumselfsimilar}, we
can assume that $E$ is the attractor of an i.f.s. of the form $\{
\zeta R_\theta x + d_i\}_{i=1}^n$, where $\theta/\pi$ is
irrational, and the open set condition is satisfied. Write
$f_i(x)= \zeta R_\theta x + d_i$.

We can assume without loss of generality that $f_i(B)\subset B$
for all $i$, where $B$ is the unit ball. Fix some small $\e>0$ and
some large $m$. We can apply Proposition
\ref{prop:discretemarstrand} to the family $\Q=\{ f_u(B):|u|=m\}$,
with $\rho = \zeta^m$, $\gamma=\dim(E)$, and some $A, A_1, A_2$
which depend only on $E$. In this way we obtain a set of ``good''
angles $J$ verifying properties (1)-(3).

Let $\rot$ be rotation by $m\theta$ on the circle. Because of the
irrationality assumption, this is a uniquely ergodic
transformation.

 Fix some $\xi\in [0,\pi)$; we
will use an inductive construction to show that
\[
\dim\left(P_\xi(E)\right) \ge \dim(E)-\eta,
\]
where $\eta$ depends on $m$ and $\e$ and can be made arbitrarily
small. Our construction will yield a tree $\R$ with vertices
labelled by disks, such that if $\R_j$ is the set of vertices at
level $j$, then the following conditions hold:
\begin{enumerate}
\item[(A)] If $D$ is the parent of $D'$, then $D'\subset D$.
\item[(B)] If $D\in \R_j$, then $D=f_u(B)$ for some word $u$ of length $m j$. In particular, $D$ has radius $\zeta^{jm}$.
\item[(C)] The family $\{ P_{\xi}(D) : D\in \R_j\}$ is disjoint and
$\zeta^{jm}$-separated.
\item[(D)] All vertices of level $j$ have the same number of children $C_j$. Moreover,
\begin{eqnarray*}
\rot^j(\xi) \notin J & \Longrightarrow & C_j =1,\\
\rot^j(\xi) \in J & \Longrightarrow & C_j \ge \e \delta
\zeta^{-m\gamma}.
\end{eqnarray*}
\end{enumerate}

Let $B$ be the root of the tree; i.e., the only vertex of level
$0$. Now assume that a vertex $D$ of level $j$ has been defined;
by property (E) we have $D=f_u(B)$, where $u$ is a word of length
$jm$. We consider two cases: if $\rot^j(\xi)\notin J$, then $D$
has just one child, $f_{u 1^m}(B)$. If $\rot^j(\xi)\in J$, then
Proposition \ref{prop:discretemarstrand} guarantees that the
family $\Q$ has a subcollection $\{ f_{v_i}(B)\}_{i=1}^M$ such
that all the projections $P_{\rot^j(\xi)}(f_{v_i}(B))$ are
$\zeta^m$-separated, and $M\ge \e \delta \zeta^{-m\gamma}$. We
define the  set of offspring of $D$ to be $\{ f_u \circ f_{v_i} (B)
\}_{i=1}^M$.

Properties (A), (B), and (D) are immediate from the construction.
To establish (C), it is enough to show that the projections of the
offspring of a given $D\in\R_j$ are $\zeta^{(j+1)m}$-separated. Let
$D=f_u(B)$; since $|u|=jm$ it follows that
\[
f_u(x) = \zeta^{jm} R_{\rot^j(0)}(x) + d_u,
\]
for some translation $d_u$. Using this we get
\begin{eqnarray*}
P_\xi(f_u \circ f_{v_i}(B)) & = & P_\xi\left(\zeta^{jm}
R_{\rot^j(0)} f_{v_i}(B) + d_u\right) \\
& = & \zeta^{jm} P_{\rot^j(\xi)}\left( f_{v_i}(B)\right) +
P_\xi(d_u).
\end{eqnarray*}
This shows that (C) holds as well.

Now let
\[
E_\xi = \bigcap_{j=1}^\infty \bigcup \{ D: D\in\R_j\}.
\]
By properties (A) and (B) this set is well defined and contained
in $P_\xi(E)$. Let $\mu_\xi$ be the probability measure which
assigns equal mass to all intervals $P_\xi(D)$ for $D\in\R_j$;
this is well defined because of properties (C) and (D); moreover,
$\mu_\xi$ is supported on $E_\xi$. We estimate the dimension of
$E_\xi$ using $\mu_\xi$ as a Frostman measure, applying Weyl's
equidistribution principle to bound from above the
$\mu_\xi$-measure of projections $P_\xi(D)$. This is done exactly
as in the proof of Theorem \ref{th:sumselfsimilar}, so we skip the
details. In the end we obtain the estimate
\[
\dim\left(P_\xi(E)\right) \ge \frac{(1-\e)(\log(\e
\delta)+m\gamma\log(1/\zeta) )}{m\log(1/\zeta)}.
\]
The right-hand side can be made arbitrarily close to $\gamma$, and
this completes the proof. \qed

\bigskip

\textit{Proof of Proposition \ref{prop:approxhomogeneous}}. We
will present the proof in the case $p=2$; the case $p=1$ is
analogous but simpler. The linear maps $A_i$ can be decomposed as
\[
A_i(x) = r_i R_{\theta_i} O_i(x),
\]
where $|r_i|<1$, $R_{\theta_i}$ is rotation by angle $\theta_i$,
and $O_i$ is either the identity or a reflection. Without loss of
generality we can assume that $E$ does not contain reflections;
i.o. all the $O_i$ are the identity. Indeed, if this is not the
case, then, say, $O_1$ is a reflection. Then we can iterate the
i.f.s. a large number of times, and then compose each of the
resulting maps which contains a reflection with $A_1$. In this way
we obtain a new self-similar set without reflections, which is
contained in $E$ and with dimension arbitrarily close to the
dimension of $E$.

Now let $\gamma$ be the only real number such that $\sum_i
r_i^\gamma = 1$. Since $E$ verifies the open set condition, this
is also the Hausdorff dimension of $E$.

Let $\{e_1,\ldots, e_n\}$ be an orthonormal basis of
$\mathbb{R}^n$. Consider the random walk which starts at $0$ and
moves from $x$ to $x+e_i$ with probability $r_i^\gamma$. Let
$\mathbb{X}_k$ be the position of this random walk after $k$
steps. The mean of $\mathbb{X}_k$ is given by
\[
\mathbb{E}\mathbb{X}_k = \sum_{i=1}^n k r_i^\gamma e_i.
\]
Let
\[
v = \sum_{i=1}^n \left(\lceil k r_i^\gamma \rceil \right) e_i,
\]
and notice that we have $|v - \mathbb{E}\mathbb{X}_k| \le
\sqrt{n}$. It follows from \cite[Chapter II, Proposition
P9]{spitzer} that
\begin{equation} \label{eq:randomwalkprob}
P(\mathbb{X}_k=v) \ge c \left(\sqrt{k}\right)^{1-n},
\end{equation}
for some $c>0$ independent of $k$. Each of the paths which end in
$v$ has probability
\[
\prod_{i=1}^n r_i^{\gamma v_i} \le  \prod_{i=1}^n r_i^{\gamma k
r_i^\gamma}.
\]
If $N_k$ is the number of such paths it follows from
(\ref{eq:randomwalkprob}) that
\begin{equation} \label{eq:estimatenumber}
N_k \ge c \left(\sqrt{k}\right)^{1-n} \prod_{i=1}^n r_i^{-k \gamma
r_i^\gamma}.
\end{equation}
Let
\begin{equation} \label{eq:estimaterho}
\rho = \prod_{i=1}^n r_i^{v_i} \ge \prod_{i=1}^n r_i \prod_{i=1}^n
r_i^{k r_i^\gamma}.
\end{equation}
We can identify paths on the lattice with compositions of maps in
the original i.f.s. Notice that all paths ending in $v$ are
associated to a map of the form $A x + z_i$, where the similarity
ratio of $A$ is $\rho$. The attractor $H$ of the i.f.s. generated
by such maps is contained in $E$, and its dimension $\tau$
verifies
\begin{equation} \label{eq:dimensionH}
N_k \rho^\tau = 1 \,\Longrightarrow \tau =
\frac{\log(N_k)}{\log(1/\rho)}.
\end{equation}
From (\ref{eq:estimatenumber}), (\ref{eq:estimaterho}) and
(\ref{eq:dimensionH}) we obtain
\[
\tau \ge \frac{\log c + (1-n)\log \sqrt{k} + k \gamma \sum_{i=1}^n
 r_i^\gamma \log(1/r_i) }{ \sum_{i=1}^n \log(1/r_i) + k \sum_{i=1}^n r_i^\gamma \log(1/r_i) }.
\]
Since the right-hand terms dominate on both numerator and
denominator, it follows that $\tau$ can be made arbitrarily close
to $\gamma$ by taking $k$ large enough, completing the proof. \qed

\section{Concluding remarks and open questions}
\label{sec:remarks}

\begin{enumerate}
\item \textbf{Sums of graph-directed attractors}.
Theorem \ref{th:sumselfsimilar} can be further generalized. The
sets $K$ and $K'$ can be attractors of more general graph-directed
systems. In order to see this, one can either observe that the
proof extends to this generality with minor variations, or just
use the fact that attractors of graph-directed systems contain
self-similar sets of arbitrarily close dimension.
\item \textbf{Resonance as an equivalence relation}. Let
$\mathcal{I}$ denote the set of all self-similar i.f.s.\ satisfying
the open set condition and such that the attractor has dimension
at most $1/2$. Algebraic resonance is clearly an equivalence
relation on $\mathcal{I}$; it follows from our results that so is
geometrical resonance. For more general compact sets there is no
clear way to define resonance algebraically, but we can say that
two compact sets $K, K'$ of Hausdorff dimension at most $1/2$ are
geometrically resonant if there is $s>0$ such that
\[
\dim(K+s K') < \dim(K)+\dim(K').
\]
However, it turns out that this is not an equivalence relation. To
see this, fix two very close numbers $a<b$ such that $C_a$ and
$C_b$ have dimension less than $1/2$; moreover, from the proof of
Theorem \ref{th:exceptions} it follows that we can ensure that the
dimension drop in $C_a+C_a$ and $C_b+C_b$ is at least $4\eta$,
where $\eta=\dim(C_b)-\dim(C_a)$. Construct a compact set $K$ in
the following way: pick a rapidly increasing sequence $n_j$.
Follow the construction of $C_a$ for $n_1$ steps, then the
construction of $C_b$ for $n_2$ steps, and so on. If this is done
carefully, then one can ensure that $K$ has dimension $\dim(C_a)$
and
\[
\max(\dim(K + C_a), \dim(K +C_b)) < 2\dim(C_b) - 2\eta,
\]
whence $K$ is resonant to both $C_a$ and $C_b$, yet $C_a$ and
$C_b$ need not be resonant.
\item \textbf{Uncountably many resonances}. We can rephrase (\ref{eq:sumcantorcompactdim})
in Theorem \ref{th:peressolomyak} as saying that a compact set
cannot be resonant to a positive measure set of central Cantor
sets. In light of the results of this paper, a natural question is
whether an arbitrary compact set is resonant to at most countably
many of the $C_a$. The answer, however, turns out to be negative.
We sketch the construction of a compact set resonant to
uncountably many central Cantor sets.

First of all, let us say that two compact sets $K, K'$ are
$\gamma$-resonant at scale $\rho$ if
\[
\leb(K_\rho + K'_\rho) \le \rho^{-\gamma},
\]
where $K_\rho, K'_\rho$ denote the $\rho$-neighborhoods of $K,
K'$. Pick $a_1>a_2$ very close to each other so that, for some
$\eta>0$,
\begin{enumerate}
\item[(i)] $\dim(C_{a_i})<1/2$.
\item[(ii)] $\dim(C_{a_i}) \in B_{\e/10}(\eta)$.
\item[(iii)]
\[
2\dim(C_{a_i}) - \dim(C_{a_i}+C_{a_i}) > 10 \e, \quad i=1,2.
\]
\end{enumerate}

We construct the desired compact set $K$ as follows (this is
different from the set constructed in the previous remark):
starting from the unit interval, follow the construction of
$C_{a_1}$ for $n_1$ steps, so that at scale $\rho_1=a_1^{n_1}$ the
sets $K$ and $C_{a_1}$ are $(2\eta-2\e)$-resonant. Observe that
there exists a small interval $I_1$ around $a_1$ such that $K$ and
$C_a$ are $(2\eta-\e)$-resonant at scale $\rho_1$ for all $a\in
I_1$; we can assume that $a_2\notin I_1$. Pick any $a_{11} <
a_{12} < a_1$ in $I_1$.

Next we follow the construction of $C_{a_2}$ for $n_2$ steps, with
$n_2>>n_1$. In this way we can find a $\rho_2<<\rho_1$ such that
$C_{a_2}$ and $K$ are $(2\eta-\e)$-resonant at scale $\rho_2$. We
find an interval $I_2$ around $a_2$, disjoint from $I_1$, and
points $a_{21}>a_{22}>a_2$ as above.

We repeat this process for $a_{ij}$ to find intervals $I_{ij}$ and
points $a_{ijk}\in I_{ij}$ such that at some very small scales,
$\rho_{ij}$ the sets $K$ and $C_{a_{ij}}$ are
$(2\eta-\e)$-resonant. We continue this construction inductively.
Then if we let
\[
\mathcal{A} = \bigcap_{k=1}^\infty\, \bigcup_{i_1,\ldots,
i_k\in\{1,2\}} I_{i_1 \ldots i_k},
\]
it follows from the construction that if $a\in \mathcal{A}$, then
$K$ and $C_a$ are $(2\eta-\e)$-resonant at arbitrarily small
scales. This implies that $C_a+K$ has dimension at most
$2\eta-\e$. Also, $\dim(K) \ge \eta-\e/10$. We conclude that
\begin{eqnarray*}
\dim(C_a+K) & < & (\eta-\e/10)+(\eta-\e/10) - 4\e/5\\
 & < & \dim(C_a)+\dim(K) -\e/2.
\end{eqnarray*}

The uncountable set of resonances in this example is of dimension
zero. We do not know if there exists a compact set $K$ which is geometrically
resonant to the central Cantor set $C_a$ for a set of parameters $a$ of positive dimension.
\item \textbf{The case $\dim(C_a)+\dim(C_b)>1$}. In this paper we
focused on the case where the sum of the dimensions is at most
$1$. We do prove that if $\dim(C_a)+\dim(C_b)>1$ and $\log b/\log
a$ is irrational, then $\dim(C_a+C_b)=1$. However, in light of
Marstrand's projection theorem and Theorem
\ref{th:peressolomyak}(ii), it is natural to conjecture that
$C_a+C_b$ actually has positive measure.

By analogy with the general results on intersections with lines
(See \cite[Chapter 10]{mattila}) and Theorem 1.2 in
\cite{peres-solomyak-cantor}, we also conjecture that, under the
irrationality assumption, for all $\theta\in (0,\pi)\backslash
\{\pi/2\}$ there is a set of positive measure of lines with
direction $\theta$ which intersect $C_a\times C_b$ in a set of
dimension $\dim(C_a)+\dim(C_b)-1$. Note that this would imply that
$C_a+C_b$ has positive measure.
\end{enumerate}

\noindent{\bf Acknowledgement}. We are grateful to Hillel
Furstenberg and Boris Solomyak for introducing us to these
problems and for many enlightening discussions. We also thank
Kemal Ilgar Ero\v{g}lu and the referee for helpful comments.


\begin{thebibliography}{10}

\bibitem{CabrelliHareMolter}
C.~A. Cabrelli, K.~E. Hare, and U.~M. Molter.
\newblock Sums of {C}antor sets yielding an interval.
\newblock {\em J. Aust. Math. Soc.}, 73(3):405--418, 2002.

\bibitem{eroglu-projections}
K.~I. Ero\v{g}lu.
\newblock On planar self-similar sets with a dense set of rotations.
\newblock {\em Annales Academiae Scientiarum Fennicae Mathematica, to appear.
  {\rm available at arXiv:math.CA/0603181}}, 2007.

\bibitem{eroglu-sums}
K.~I. Ero\v{g}lu.
\newblock On the arithmetic sums of cantor sets.
\newblock {\em Nonlinearity}, 20:1145--1161, 2007.

\bibitem{falconer2}
K.~Falconer.
\newblock {\em Fractal geometry}.
\newblock John Wiley \& Sons Ltd., Chichester, 1990.
\newblock Mathematical foundations and applications.

\bibitem{falconer3}
K.~Falconer.
\newblock {\em Techniques in fractal geometry}.
\newblock John Wiley \& Sons Ltd., Chichester, 1997.

\bibitem{furstenberg}
H.~Furstenberg.
\newblock Intersections of {C}antor sets and transversality of semigroups.
\newblock In {\em Problems in analysis (Sympos. Salomon Bochner, Princeton
  Univ., Princeton, N.J., 1969)}, pages 41--59. Princeton Univ. Press,
  Princeton, N.J., 1970.

\bibitem{keane-simon-solomyak}
M.~Keane, K.~Simon, and B.~Solomyak.
\newblock The dimension of graph directed attractors with overlaps on the line,
  with an application to a problem in fractal image recognition.
\newblock {\em Fund. Math.}, 180(3):279--292, 2003.

\bibitem{keane-smorodinsky-solomyak}
M.~Keane, M.~Smorodinsky, and B.~Solomyak.
\newblock On the morphology of {$\gamma$}-expansions with deleted digits.
\newblock {\em Trans. Amer. Math. Soc.}, 347(3):955--966, 1995.

\bibitem{kemeny-snell}
J.~G. Kemeny and J.~L. Snell.
\newblock {\em Finite {M}arkov chains}.
\newblock The University Series in Undergraduate Mathematics. D. Van Nostrand
  Co., Inc., Princeton, N.J.-Toronto-London-New York, 1960.

\bibitem{kenyon-projecting}
R.~Kenyon.
\newblock Projecting the one-dimensional {S}ierpinski gasket.
\newblock {\em Israel J. Math.}, 97:221--238, 1997.

\bibitem{lagarias-wang}
J.~C. Lagarias and Y.~Wang.
\newblock Tiling the line with translates of one tile.
\newblock {\em Invent. Math.}, 124(1-3):341--365, 1996.

\bibitem{mattila}
P.~Mattila.
\newblock {\em Geometry of sets and measures in {E}uclidean spaces}, volume~44
  of {\em Cambridge Studies in Advanced Mathematics}.
\newblock Cambridge University Press, Cambridge, 1995.
\newblock Fractals and rectifiability.

\bibitem{MendesOliveira}
P.~Mendes and F.~Oliveira.
\newblock On the topological structure of the arithmetic sum of two {C}antor
  sets.
\newblock {\em Nonlinearity}, 7(2):329--343, 1994.

\bibitem{moreira-hungarica}
C.~G. T. d.~A. Moreira.
\newblock Sums of regular {C}antor sets, dynamics and applications to number
  theory.
\newblock {\em Period. Math. Hungar.}, 37(1-3):55--63, 1998.
\newblock International Conference on Dimension and Dynamics (Miskolc, 1998).

\bibitem{moreirayoccoz}
C.~G. T. d.~A. Moreira and J.-C. Yoccoz.
\newblock Stable intersections of regular {C}antor sets with large {H}ausdorff
  dimensions.
\newblock {\em Ann. of Math. (2)}, 154(1):45--96, 2001.

\bibitem{palis-conjecture}
J.~Palis.
\newblock Homoclinic orbits, hyperbolic dynamics and dimension of {C}antor
  sets.
\newblock In {\em The Lefschetz centennial conference, Part III (Mexico City,
  1984)}, volume~58 of {\em Contemp. Math.}, pages 203--216. Amer. Math. Soc.,
  Providence, RI, 1987.

\bibitem{peresschlag}
Y.~Peres and W.~Schlag.
\newblock Smoothness of projections, {B}ernoulli convolutions, and the
  dimension of exceptions.
\newblock {\em Duke Math. J.}, 102(2):193--251, 2000.

\bibitem{sixtyyears}
Y.~Peres, W.~Schlag, and B.~Solomyak.
\newblock Sixty years of {B}ernoulli convolutions.
\newblock In {\em Fractal geometry and stochastics, II (Greifswald/Koserow,
  1998)}, volume~46 of {\em Progr. Probab.}, pages 39--65. Birkh\"auser, Basel,
  2000.

\bibitem{peres-solomyak-cantor}
Y.~Peres and B.~Solomyak.
\newblock Self-similar measures and intersections of {C}antor sets.
\newblock {\em Trans. Amer. Math. Soc.}, 350(10):4065--4087, 1998.

\bibitem{pollicott-simon}
M.~Pollicott and K.~Simon.
\newblock The {H}ausdorff dimension of {$\lambda$}-expansions with deleted
  digits.
\newblock {\em Trans. Amer. Math. Soc.}, 347(3):967--983, 1995.

\bibitem{rams-pacific}
M.~Rams.
\newblock Generic behavior of iterated function systems with overlaps.
\newblock {\em Pacific J. Math.}, 218(1):173--186, 2005.

\bibitem{shmerkin-moreira}
P.~Shmerkin.
\newblock {M}oreira's theorem on the sum of regular {C}antor sets.
\newblock {\em Preprint}, 2006.

\bibitem{solomyak-indagationes}
B.~Solomyak.
\newblock On the measure of arithmetic sums of {C}antor sets.
\newblock {\em Indag. Math. (N.S.)}, 8(1):133--141, 1997.

\bibitem{spitzer}
F.~Spitzer.
\newblock {\em Principles of random walk}.
\newblock Springer-Verlag, New York, second edition, 1976.
\newblock Graduate Texts in Mathematics, Vol. 34.

\end{thebibliography}
\end{document}